\documentclass[a4paper,12pt]{article}
\usepackage{graphics}
\usepackage{amsmath}
\usepackage{cases}
\usepackage{mathrsfs}
\usepackage{amsfonts}
\usepackage{latexsym}
\usepackage{indentfirst}

\makeatletter

\newcommand{\Rmnum}[1]{\expandafter\@slowromancap\romannumeral #1@}
\title{Bicyclic graphs with exactly two main signless Laplacian eigenvalues
\thanks{This work was supported by Hunan Provincial Natural Science Foundation of China.}}
\author{{\sc He Huang, Hanyuan Deng}\thanks{Corresponding author:
hydeng@hunnu.edu.cn.}\\
{\small College of Mathematics and Computer Science,}\\
{\small Hunan Normal University, Changsha, Hunan 410081, P. R.
China}}
\date{}
\begin{document}
\maketitle

\begin{abstract}
A signless Laplacian eigenvalue of a graph $G$ is called a main
signless Laplacian eigenvalue if it has an eigenvector the sum of
whose entries is not equal to zero. In this paper, all connected
bicyclic graphs with exactly two main eigenvalues are determined.



\end{abstract}

\baselineskip=0.30in

\section{Introduction}
Let $G$ be a simple connected graph with vertex set $V=\{v_1, v_2$,
$\cdots, v_n\}$ and $(0,1)$-adjacency matrix $A$. An eigenvalue
$\lambda$ of $A$ is said to be a main eigenvalue of $G$ if the
eigenspace $\varepsilon(\lambda)$ is not orthogonal to the all-1
vector $\mathbf{j}$, i.e., it has an eigenvector the sum of whose
entries is not equal to zero. An eigenvector $\mathbf{x}$ is a main
eigenvector if $\mathbf{x}^{T}\mathbf{j}\neq 0$. For details of the
main eigenvector of adjacency matrix, the readers are suggested to
refer to \cite{cv,ha,pr}.

Very recently, we introduced in~\cite{d} the main signless Laplacian
eigenvalue of a graph. Let $L^{+}=D+A$, where
$D=diag(d(v_1),d(v_2),\cdots,d(v_n))$ is the diagonal matrix of
vertex degrees and $d(v_i)$ is the degree of vertex $v_i$. Then
$L^+$ is called the signless Laplacian matrix of $G$. Some results
on the signless Laplacian matrix of a graph can be found
in~\cite{dc4,dc1,dc5,dc2,dc3}. An eigenvalue $\mu$ of $L^{+}$ is
said to be a main signless Laplacian eigenvalue of $G$ (or a main
eigenvalue of $L^+$) if the eigenspace $\varepsilon(\mu)$ is not
orthogonal to the all-1 vector $\mathbf{j}$, i.e., it has an
eigenvector the sum of whose entries is not equal to zero. It was
showed in~\cite{d} that a graph $G$ with exactly one main signless
Laplacian eigenvalue if and only if $G$ is regular, trees and
unicyclic graphs with exactly two main signless Laplacian
eigenvalues were characterized.

A graph $G$ is called a bicyclic graph, if it is a simple connected
graph in which the number of edges equals the number of vertices
plus one. The aim of this work is to characterize all bicyclic
graphs with exactly two main signless Laplacian eigenvalues.

\section{Some Lemmas}
In~\cite{d}, an alternative characterization of a graph with exactly
two main signless Laplacian eigenvalues is given by 2-walks. A graph
$G$ is 2-walk {\bf $(a, b)$-parabolic}~\cite{d} if there exist
uniquely a positive integer $a$ and a non-negative integer $b$ such
that such that $a^2-8b>0$ and $s(v)=-d^2(v)+ad(v)-b$ holds for all
$v\in V(G)$, where $s(v)=\sum_{u\in N_G(v)}d(u)$ and $N_G(v)$ is the
set of all neighbors of $v$ in $G$. Note that any 2-walk $(a,
b)$-parabolic graph is irregular, since one has $s(v)
=-d^2(v)+2kd(v)-0$ and $s(v)=-d^2(v)+(2k+1)d(v)-k$ for a $k$-regular
graph, i.e., $(a,b)=(2k,0)$ or $(2k+1, k)$ is not unique.

{\bf Lemma 1}~\cite{d}. A graph $G$ has exactly two main signless
Laplacian eigenvalues if and only if there exist uniquely a positive
integer $a$ and a non-negative integer $b$ such that $a^2-8b>0$ and
\begin{equation}\label{eq:1}
s(v)=-d^2(v)+ad(v)-b
\end{equation}
for all $v\in V(G)$, i.e., $G$ is a 2-walk $(a, b)$-parabolic graph.

In~\cite{d}, we showed that a graph has exactly one main signless
Laplacian eigenvalue if and only if it is regular. So, a 2-walk $(a,
b)$-parabolic graph $G$ is non-regular, and there are $u, v\in V(G)$
such that $d(u)\neq d(v)$. From $s(u)=-d^2(u)+ad(u)-b$ and
$s(v)=-d^2(v)+ad(v)-b$, we have
\begin{equation}\label{eq:2}
a=\frac{s(u)-s(v)}{d(u)-d(v)}+d(u)+d(v),\,\,\,\,\,\,\,
b=\frac{s(u)-s(v)}{d(u)-d(v)}d(v)+d(u)d(v)-s(v)
\end{equation}

{\bf Remark 1}. If $G$ is a 2-walk $(a, b)$-parabolic graph with
$\delta(G)=1$, then $a-b\geq 3$ since there is a pendent vertex $x$
with the only incident edge $xy$ in $G$ and $d(y)=s(x)=-1+a-b\geq
2$.

Let $\mathscr{G}_{a,b}=\{G: G \mbox{ is a 2-walk } (a,
b)-\mbox{parabolic bicyclic graph with}$ $\delta(G)=1\}$. For each
$G\in \mathcal{G}_{a,b}$, $G_0$ is the graph obtained from $G$ by
deleting all pendant vertices. If $v\in V(G_0)$, we use $d_{G_0}(v)$
to denote the degree of the vertex $v$ in $G_0$. The following
lemmas are true for unicyclic graphs, and it is easy to see that
they are also true for bicyclic graphs from their proofs
in~\cite{d}.

{\bf Lemma 2}~\cite{d}. If $G\in\mathscr{G}_{a,b}$ and $v\in
V(G_0)$, then $d(v)=d_{G_0}(v)$ or $d(v)=a-b-1$.

{\bf Lemma 3}~\cite{d}. If $G\in\mathscr{G}_{a, b}$, then (i)
$\delta(G_0)\geq 2$; (ii) $a-b\geq 4$ and $a\geq 5$.

{\bf Lemma 4}. Let $G\in \mathscr{G}_{a,b}$ and $C=u_1u_2\ldots
u_ku_1$ a cycle of $G$. If $d_{G_0}(u_1)\geq3$ and $d_{G_0}(u_2)=2$,
then there is $i\in\{1,2,\ldots,k\}$ such that $d(u_i)\neq a-b-1$.

{\bf Proof}. By the definition of $G_0$, $C=u_1u_2\ldots u_ku_1$ is
also a cycle of $G_0$.

If $d(u_1)=d(u_2)=\cdots=d(u_k)=a-b-1$, then from (\ref{eq:1})
\begin{equation}\label{eq:3}
\qquad s(u_1)=-(a-b-1)^2+a(a-b-1)-b=s(u_2)
\end{equation}
And
\begin{align}
\notag s(u_1)&=d(u_1)-d_{G_0}(u_1)+d(u_2)+d(u_k)+\sum_{w\in
N_{G_0}(u_1)\backslash \{u_2,u_k\}}d(w)\\
\notag &\geq3(a-b-1)+2(d_{G_0}(u_1)-2)-d_{G_0}(u_1)=3(a-b-1)+d_{G_0}(u_1)-4\\
\notag &\geq 3(a-b-1)-1 \qquad  (\mbox{by } d_{G_0}(u_1)\geq3)\\
\notag s(u_2)&=d(u_1)+d(u_3)+d(u_2)-d_{G_0}(u_2) \qquad (\mbox{by } d_{G_0}(u_2)=2) \\
\notag       &=3(a-b-1)-2
\end{align}
So, $s(u_1)>s(u_2)$, it contradicts (\ref{eq:3}). \hfill $\Box$

If $R$ is a cycle or a path of $G$, then the length of $R$, denoted
by $l(R)$, is defined as the number of edges of $R$.

{\bf Lemma 5}. Let $G$ be a 2-walk $(a, b)$-parabolic graph.
$R=x_1x_2\cdots x_t$ is a path or cycle of length at least 2 in $G$
such that $d(x_1)\geq3$, $d(x_t)\geq3$ and
$d(x_2)=\cdots=d(x_{t-1})=2$. Then (i) If $d(x_1)=d(x_t)$, then
$l(R)\leq3$; Moreover, if $l(R)=3$, then there exists no path
$Q=y_1y_2y_3$ in $G$ such that $d(y_1)=d(y_3)=d(x_1)$ and
$d(y_2)=2$; (ii) If $d(x_1)\neq d(x_t)$, then $l(R)\leq2$.

{\bf Proof}. (i) If $l(R)\geq4$, then $d(x_2)=d(x_3)=d(x_4)=2$. By
(\ref{eq:1}), we have
$s(x_2)=d(x_3)+d(x_1)=-2^2+2a-b=s(x_3)=d(x_2)+d(x_4)=4$, and
$d(x_1)=2$. This contradicts $d(x_1)\geq3$. So, $l(R)\leq3$.

If $l(R)=3$, and there is a path $Q=y_1y_2y_3$ in $G$ such that
$d(y_1)=d(y_3)=d(x_1)$ and $d(y_2)=2$, then by (\ref{eq:1}),
$s(y_2)=d(y_1)+d(y_3)=-2^2+2a-b=s(x_2)=d(x_1)+d(x_3)=d(x_1)+2$, and
$d(y_1)+d(y_3)=d(x_1)+2$. This is impossible since
$d(y_1)=d(y_3)=d(x_1)\geq 3$.

(ii) If $l(R)\geq3$, then $s(x_2)=-2^2+2a-b=s(x_{t-1})$ by
(\ref{eq:1}). On the other hand, $s(x_2)=d(x_1)+d(x_3)=d(x_1)+2$ and
$s(x_{t-1})=d(x_{t-2})+d(x_t)=d(x_t)+2$, a contradiction since
$d(x_1)\neq d(x_t)$. \hfill $\Box$

{\bf Lemma 6}. Let $G\in\mathscr{G}_{a,b}$, $R=u_1u_2\cdots u_k$ a
path or cycle of $G_0$ with length at least 3 such that
$d_{G_0}(u_1)=d_{G_0}(u_k)\in\{3,4\}$ and $d_{G_0}(u_i)=2$ for
$2\leq i\leq k-1$.

(i) If there is $i\in\{2,\ldots,k-1\}$ such that $d(u_i)\neq
d(u_{i+1})$, then $l(R)=3$ and $b=1$. Moreover, if $d_{G_0}(u_1)=3$,
then $a=6$, $d(u_1)=d_{G_0}(u_1)=3$, $d(u_2)=2$, $d(u_3)=4$; or
$d(u_2)=4$, $d(u_3)=2$. If $d_{G_0}(u_1)=4$, then $a=7$,
$d(u_1)=d_{G_0}(u_1)=4$, $d(u_2)=2$, $d(u_3)=5$; or $d(u_2)=5$,
$d(u_3)=2$.

(ii) If $d(u_2)=d(u_3)=\cdots=d(u_{k-1})=d$, then $d\in\{2,a-b-1\}$.
Moreover, if $R$ is a cycle, then $l(R)=3$ and $d(u_2)=d(u_3)=2$.

{\bf Proof}. (i) By Lemma 2,
$d(u_j)\in\{d_{G_0}(u_j),a-b-1\}=\{2,a-b-1\}$, $2\leq j\leq k-1$. If
there is $i\in\{2,3,\cdots,k-2\}$ such that $d(u_i)\neq d(u_{i+1})$,
then $d(u_i)=2$ and $d(u_{i+1})=a-b-1$, or $d(u_{i+1})=2$ and
$d(u_i)=a-b-1$. Without loss of generality, assume that $d(u_i)=2$,
$d(u_{i+1})=a-b-1$. From (\ref{eq:2}), we have
\begin{align}
\notag
a=\frac{s(u_{i+1})-s(x)}{d(u_{i+1})-d(x)}+d(u_{i+1})+d(x)=\frac{2+d(u_{i+2})+a-b-3-(a-b-1)}{a-b-2}+a-b
\end{align}
where $x$ is a pendent vertex of $G$, and
\begin{equation}\label{eq:4}
d(u_{i+2})=b(a-b-2)
\end{equation}
$b\neq 0$ since $d(u_{i+2})\neq 0$. So, $b\geq1$.

{\bf Case 1}. $b=1$. Then $d(u_{i+2})=a-3<a-2=a-b-1$, and
$d(u_{i+2})=d_{G_0}(u_{i+2})$ from Lemma 2. If $u_{i+2}\neq u_k$,
then $d(u_{i+2})=d_{G_0}(u_{i+2})=2$. So, $a=5$, $a-b-1=3$. Because
$d(u_1)\geq d_{G_0}(u_1)\in\{3,4\}$ and $d(u_1)\in\{d_{G_0}(u_1),
a-b-1\}=\{d_{G_0}(u_1),3\}$, by Lemma 2, we have
$d(u_1)=d_{G_0}(u_1)$ and $N_{G}(u_1)=N_{G_0}(u_1)$.

From (\ref{eq:1}), $s(u_1)=-d^2_{G_0}(u_1)+5d_{G_0}(u_1)-1$. On the
other hand, $s(u_1)=\sum\limits_{w\in
N_G(u_1)}d(w)=\sum\limits_{w\in N_{G_0}(u_1)}d(w)\geq
2d_{G_0}(u_1)$. So,
$-d^2_{G_0}(u_1)+5d_{G_0}(u_1)-1\geq2d_{G_0}(u_1)$. This implies
$\frac{3-\sqrt{5}}{2}\leq d_{G_0}(u_1)\leq \frac{3+\sqrt{5}}{2}$, a
contradiction with $d_{G_0}(u_1)\geq3$. So, $u_{i+2}=u_k$, and
$k=i+2$.

Since $d(u_i)=2$, $s(u_i)=-2^2+2a-b=2a-5$ from (\ref{eq:1}), and
$d(u_{i-1})=s(u_i)-d(u_{i+1})=2a-5-(a-2)=a-3\neq a-b-1$. By Lemma 2,
$d_{G_0}(u_{i-1})=d(u_{i-1})=a-3=d(u_{i+2})=d(u_k)>2$. So,
$u_{i-1}=u_1$, and $i=2$, $l(R)=k-1=i+2=3$.

If $d_{G_0}(u_1)=3$, then
$d_{G_0}(u_1)=d_{G_0}(u_{i-1})=d(u_{i-1})=a-3$ from above, and
$a=6$. $d(u_2)=d(u_i)=2$ and $d(u_3)=d(u_{i+1})=a-b-1=4$.

Similarly, if $d_{G_0}(u_1)=4$, then $a-3=4$, and $a=7$. $d(u_2)=2$
and $d(u_3)=5$.

{\bf Case 2}. $b=2$.

If $d_{G_0}(u_k)=3$, then from (\ref{eq:4}) and Lemma 3,
$d(u_{i+2})=2(a-b-2)=a-b-1+a-b-3>a-b-1\geq 3\geq d_{G_0}(u_{i+2})$.
This contradicts Lemma 2 whether $R=u_1u_2\cdots u_k$ a path or
cycle. If $d_{G_0}(u_k)=4$, then from (\ref{eq:4}) and Lemma 3,
$d(u_{i+2})=2(a-b-2)>a-b-1>2$. By Lemma 2, we have
$d_{G_0}(u_{i+2})=d(u_{i+2})>2$. So, $u_{i+2}=u_k$, and
$d(u_k)=d_{G_0}(u_k)=4$. This implies $2(a-b-2)=4$ and $a=6$.

From (\ref{eq:1}), $s(u_k)=-d^2(u_k)+ad(u_k)-b=6$. On the other
hand, $s(u_k)=\sum_{w\in N_G(u_k)}d(w)=\sum_{w\in
N_{G_0}(u_k)}d(w)\geq 2d_{G_0}(u_k)=8$, a contradiction.

{\bf Case 3}. $b>2$.

From (\ref{eq:4}), $d(u_{i+2})=b(a-b-2)\geq
3(a-b-2)=a-b-1+2(a-b-2)-1>max\{a-b-1, d_{G_0}(u_{i+2})\}$. This
contradicts Lemma 2.

(ii) It is obvious that $d(u_2)=\cdots=d(u_{k-1})\in\{2,a-b-1\}$
from Lemma 2, since $d(u_2)=\cdots=d(u_{k-1})$ and
$d_{G_0}(u_2)=\cdots=d_{G_0}(u_{k-1})=2$.

Moreover, if $R$ is a cycle and $d(u_2)=\cdots=d(u_{k-1})=a-b-1$,
then $d(u_1)\neq a-b-1$ from Lemma 4, and $d(u_1)=d_{G_0}(u_1)$ from
Lemma 2. By (\ref{eq:1}), $s(u_2)=-(a-b-1)^2+a(a-b-1)-b$. On the
other hand,
$s(u_2)=d(u_1)+d(u_3)+d(u_2)-d_{G_0}(u_2)=2(a-b-1)-2+d(u_1)$. So,
\begin{align}
\notag d(u_1)&=s(u_2)-2(a-b-1)+2=-(a-b-1)^2+a(a-b-1)-b-2(a-b-1)+2\\
\notag       &=(b-1)(a-b-1)-b+2=(b-1)(a-b-2)+1
\end{align}

If $b\leq1$, then $d(u_1)\leq1$. This is impossible. If $b=2$, then
$d(u_1)=a-b-1$. This contradicts $d(u_1)\neq a-b-1$. If $b>2$, then
$d(u_1)\geq 2(a-b-2)+1\geq 5>d_{G_0}(u_1)$. But $d(u_1)\neq a-b-1$
from Lemma 4, it contradicts Lemma 2.

So, $d(u_2)=\cdots=d(u_{k-1})=2$. And $l(R)=3$ from Lemma 5(i).
\hfill $\Box$

{\bf Lemma 7}. Let $G\in\mathscr{G}_{a,b}$, $u, v\in V(G_0)$ and
$d(u)=d(v)$. If $N_{G_0}(u)=\{u',u''\}$ and $N_{G_0}(v)=\{v',v''\}$,
then $d(u')=d(v')$ if and only if $d(u'')=d(v'')$.

{\bf Proof}. From (\ref{eq:1}) and $d(u)=d(v)$, we have
$s(u)=-d^2(u)+ad(u)-b=-d^2(v)+ad(v)-b=s(v)$. On the other hand,
$s(u)=d(u')+d(u'')+d(u)-d_{G_0}(u)$ and
$s(v)=d(v')+d(v'')+d(v)-d_{G_0}(v)$. So, $d(u')=d(v')$ if and only
if $d(u'')=d(v'')$. \hfill $\Box$

\begin{figure}
\includegraphics{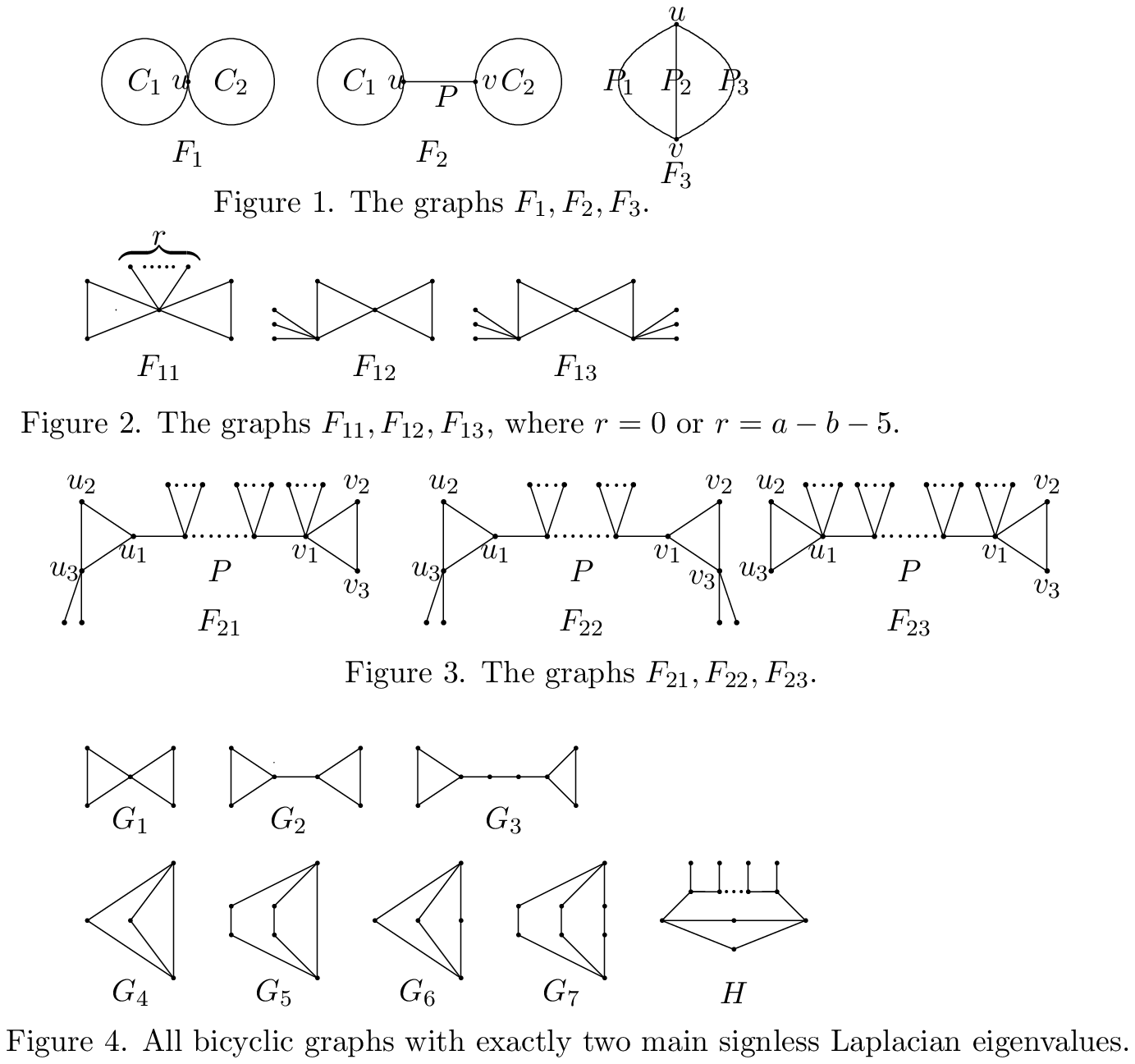}
\end{figure}

\section{Bicyclic graphs with exactly two main signless Laplacian eigenvalues}
In the following, we determine all bicyclic graphs with exactly two
main signless Laplacian eigenvalues.

{\bf Theorem 1}. A bicyclic graph $G$ with $\delta(G)=1$ has exactly
two main signless Laplacian eigenvalues if and only if $G\cong H$
(see Figure 4).

{\bf Proof}. By way of contradiction, assume that
$G\in\mathscr{G}_{a,b}$. From Lemma 3, $\delta(G_0)\geq2$. So, $G_0$
is one of the three graphs $F_1, F_2, F_3$ shown in Figure 1.

{\bf Case 1}. $G_0=F_1$, where $|V(C_1)\cap V(C_2)|=1$ ,
$C_1=u_1u_2\cdots u_n$, $C_2=v_1v_2\cdots v_m$ and
$u_1=u_n=v_1=v_m=u$.

From Lemma 6, we know that $n=m=4$ and (i) if $d(u_2)=d(u_3)$,
$d(v_2)=d(v_3)$, then $G=F_{11}$; (ii) if $d(u_2)=d(u_3)$,
$d(v_2)\neq d(v_3)$ (or $d(u_2)\neq d(u_3)$, $d(v_2)=d(v_3)$), then
$G=F_{12}$; (iii) if $d(u_2)\neq d(u_3)$, $d(v_2)\neq d(v_3)$, then
$G=F_{13}$.

Note that $G\in\mathscr{G}_{a,b}$ and $\delta(G)=1$, there is a
pendant vertex $x$ in $G$.

If $G=F_{11}$, then $d(u_2)=d(u_3)=d(v_2)=d(v_3)=2$, $d(u_1)=a-b-1$.
By (\ref{eq:2}), $$a=\frac{s(u_2)-s(x)}{d(u_2)-d(x)}+d(u_2)+d(x)=5,
a=\frac{s(u_1)-s(x)}{d(u_1)-d(x)}+d(u_1)+d(x)=\frac{4}{a-b-2}+a-b.
$$ So, $b^2-3b+4=0$. This is impossible.

If $G=F_{12}$ or $F_{13}$, then $a=7$ and $b=1$ by Lemma 6(i). Using
Eq.(\ref{eq:1}), it is easy to check that $G=F_{12}$ and $F_{13}$
are not 2-walk $(a, b)$-parabolic graphs, a contradiction.

{\bf Case 2}. $G_0=F_2$, where $V(C_1)\cap V(C_2)=\emptyset$,
$C_1=u_1u_2\cdots u_n$, $C_2=v_1v_2\cdots v_m$, $P=w_1w_2\cdots
w_t$, $u_1=u_n=w_1=u$ and $w_t=v_1=v_m=v$.

From Lemma 6, we know that $n=m=4$ and (i) if $d(u_2)\neq d(u_3)$,
$d(v_2)=d(v_3)$ (or $d(u_2)=d(u_3)$, $d(v_2)\neq d(v_3)$), then
$G\cong F_{21}$; (ii) if $d(u_2)\neq d(u_3)$, $d(v_2)\neq d(v_3)$,
then $G\cong F_{22}$; (iii) if $d(u_2)=d(u_3)$, $d(v_2)=d(v_3)$,
then $G\cong F_{23}$.

(i) In the graph $F_{21}$, we have $a=6$, $b=1$, $d(u_2)=2$,
$d(u_3)=4$, $d(u_1)=3$, $d(v_2)=2$, $d(v_3)=2$ by Lemma 6.
$d(v_1)=5$ by using Eq.(\ref{eq:1}) to $v_2$, and it can be found
that a pendant vertex $x$ of $v_1$ does not satisfy Eq.(\ref{eq:1}).
So, $F_{21}$ is not a 2-walk $(a, b)$-parabolic graph.

(ii) In the graph $F_{22}$, we have $a=6$, $b=1$, $d(u_2)=2$,
$d(u_3)=4$, $d(u_1)=3$, $d(v_2)=2$, $d(v_3)=4$, $d(v_1)=3$ by Lemma
6. First, $d(w_2)=2$ by using Eq.(\ref{eq:1}) to $u_1$, and then
$d(w_3)=4$ by using Eq.(\ref{eq:1}) to $w_2$, and finally $d(w_4)=3$
by using Eq.(\ref{eq:1}) to $w_3$. If $w_4=v_1$, then
$s(v_1)=d(w_3)+d(v_2)+d(v_3)=4+2+4\neq-d^2(v_1)+ad(v_1)-b$, i.e.,
$v_1$ does not satisfy Eq.(\ref{eq:1}); otherwise, a pendant vertex
$x$ of $w_4$ does not satisfy Eq.(\ref{eq:1}). So, $F_{22}$ is not a
2-walk $(a, b)$-parabolic graph.

(iii) In the graph $F_{23}$, $d(u_2)=d(u_3)=d(v_2)=d(v_3)=2$.

First, we show that $d(u_1)=d(v_1)=3$. If $d(u_1)\neq3$, then
$d(u_1)>d_{G_0}(u_1)=3$ and $d(u_1)=a-b-1$ from Lemma 2. By
Eq.(\ref{eq:1}), we have $s(u_2)=-2^2+2a-b=2a-b-4$. On the other
hand, $s(u_2)=d(u_1)+d(u_3)=a-b-1+2$. So, $a-b-1+2=2a-b-4$ and
$a=5$. It implies that $b=0$ since $b$ is a nonnegative integer and
$a-b-1=d(u_1)>d_{G_0}(u_1)=3$, and $d(u_1)=a-b-1=4$. Using
Eq.(\ref{eq:1}) again,
$s(u_1)=-d^2(u_1)+ad(u_1)-b=-16+20-0=4<d(u_2)+d(u_3)+d(w_2)+(d(u_1)-d_{G_0}(u_1)$,
a contradiction. So, $d(u_1)=3$. Similarly (or by symmetry),
$d(v_1)=3$.

Now, $s(u_2)=-2^2+2a-b=2a-4-b$ by Eq.(\ref{eq:1}), and
$s(u_2)=d(u_1)+d(u_3)=5$. So, $2a-b=9$. And $a=5$, $b=1$ from Lemma
3. Using Eq.(\ref{eq:1}) again, we have
$s(u_1)=5=d(u_2)+d(u_3)+d(w_2)$, and $d(w_2)=1$, a contradiction.

So, $F_{23}$ is not a 2-walk $(a, b)$-parabolic graph.

{\bf Case 3}. $G_0=F_3$, where $P_1=u_1u_2\cdots u_n$,
$P_2=v_1v_2\cdots v_m$, $P_3=w_1w_2\cdots w_t$, $u_1=v_1=w_1=u$,
$u_n=v_m=w_t=v$. And we may assume that $n,m\geq 3$ since $G$ is a
simple graph.

{\bf Claim 1}. $d(u_2)=\cdots=d(u_{n-1})$,
$d(v_2)=\cdots=d(v_{m-1})$, $d(w_2)=\cdots=d(w_{t-1})$.

If there $i\in\{2,\cdots,n-2\}$ such that $d(u_i)\neq d(u_{i+1})$,
then by Lemma 6(i), $l(P_1)=3$ (i.e., $n=4$), $b=1$, $a=6$,
$d(u_1)=d_{G_0}(u_1)=3$, $d(u_4)=d_{G_0}(u_4)=3$, $d(u_2)=2$,
$d(u_3)=4$ or $d(u_2)=4$, $d(u_3)=2$. By symmetry, we may assume
that $d(u_2)=2$, $d(u_3)=4$. Using Eq.(\ref{eq:1}),
$s(u_4)=8=d(u_3)+d(v_{m-1})+d(w_{t-1})$ implies
$d(v_{m-1})=d(w_{t-1})=2$. It can be obtained that
$d(v_{m-2})=d(w_{t-2})=4$ by using Lemma 7 for $\{u_2,v_{m-1}\}$ and
$\{u_2,w_{t-1}\}$, respectively. And $l(P_2)=l(P_3)=3$ (i.e.,
$m=t=4$) from Lemma 6(i). Now,
$s(u_1)=d(u_2)+d(v_2)+d(w_2)=10\neq-d^2(u_1)+ad(u_1)-b$, it
contradicts Eq.(\ref{eq:1}). So, $d(u_2)=\cdots=d(u_{n-1})$.

Similarly, $d(v_2)=\cdots=d(v_{m-1})$ and
$d(w_2)=\cdots=d(w_{t-1})$. \hfill $\Box$

{\bf Claim 2}. $d(u_1)=d(u_n)$.

By way of contradiction, assume that $d(u_1)\neq d(u_n)$. $d(u_1),
d(u_n)\in\{3,a-b-1\}$ from Lemma 2 and
$d_{G_0}(u_1)=d_{G_0}(u_n)=3$. By symmetry, we may assume that
$d(u_1)=a-b-1$, $d(u_n)=3$. Then $a-b-1\geq 4$, i.e.,
\begin{equation}\label{eq:5}
a-b\geq 5
\end{equation}

From Claim 1, $d(u_2)=d(u_{n-1})$, $d(v_2)=d(v_{m-1})$,
$d(w_2)=d(w_{t-1})$. And
\begin{align}
\notag s(u_1)-s(u_n)& = \left\{
\begin{array}{ll}
0, & \mbox{ $l(P_3)=1$}; \\
a-b-4, & \mbox{$l(P_3)\geq2$}.
\end{array} \right.
\end{align}

From Eq.(\ref{eq:2}), we have
\begin{align}
\qquad
 \notag a&=\frac{s(u_1)-s(u_n)}{d(u_1)-d(u_n)}+d(u_1)+d(u_n)=\left\{\begin{array}{ll}
a-b+2,& \mbox{$l(P_3)=1$};\\
a-b+3,& \mbox{$l(P_3)\geq2$.}
\end{array}\right.
\end{align}
So, $b=2$ for $l(P_3)=1$ and $b=3$ for $l(P_3)\geq2$.

From Eq.(\ref{eq:2}) again,
$s(u_1)=-(a-b-1)^2+a(a-b-1)-b=(b+1)(a-b-1)-b$. On the other hand,
$s(u_1)=(d(u_1)-d_{G_0}(u_1))+d(u_2)+d(v_2)+d(w_2)$. So,
\begin{align}
\notag
d(u_2)+d(v_2)+d(w_2)= & s(u_1)-d(u_1)+d_{G_0}(u_1)=b(a-b-1)-b+3\\
=&\left\{ \begin{array}{ll}
2(a-b-1)+1, & \mbox{ $l(P_3)=1$}; \\
3(a-b-1), & \mbox{$l(P_3)\geq2$.}
\end{array}
\right.
\end{align}

By Lemma 2, $d(u_2)\in\{d_{G_0}(u_2),a-b-1\}=\{2,a-b-1\};
d(v_2)\in\{d_{G_0}(v_2),a-b-1\}=\{2,a-b-1\};
d(w_2)\in\{d_{G_0}(w_2),a-b-1\}=\{2,3,a-b-1\}$.

If $l(P_3)\geq 2$, then Eq.(6) implies $d(u_2)=d(v_2)=d(w_2)=a-b-1$;

If $l(P_3)=1$, then $w_2=u_n$ and $d(w_2)=3$. Eq.(6) implies
$d(u_2)+d(v_2)=(a-b-1)+(a-b-3)$. So, $\{d(u_2),
d(v_2)\}=\{2,a-b-1\}$. Without loss of generality, we assume that
$d(u_2)=a-b-1$.

In the following, we show that $l(P_1)=2$. Otherwise,
$a-b-1=d(u_2)=d(u_3)=\cdots=d(u_{n-1})$ from Claim 1. By
Eq.(\ref{eq:1}), $s(u_2)=s(u_{n-1})$ since $d(u_2)=d(u_{n-1})$. But
$s(u_2)=(d(u_2)-d_{G_0}(u_2))+d(u_3)+d(u_1)$ and
$s(u_{n-1})=(d(u_2)-d_{G_0}(u_2))+d(u_{n-2})+d(u_n)$, it implies
$d(u_1)=d(u_n)$, a contradiction. So, $l(P_1)=2$, and $n=3$,
$d(u_3)=d(u_n)=3$.

Now, by (\ref{eq:1}), $s(u_2)=-(a-b-1)^2+a(a-b-1)-b=(b+1)(a-b-1)-b$.
On the other hand,
$s(u_2)=d(u_1)+d(u_3)+d(u_2)-d_{G_0}(u_2)=2(a-b-1)+1$. So,
$(b+1)(a-b-1)-b=2(a-b-1)+1$, and $(b-1)(a-b-2)=2$. This is
impossible since $b\geq 2$ and $a-b\geq 5$ from Eq.(\ref{eq:5}).
\hfill $\Box$

{\bf Claim 3}. If $d(u_1)=d(u_n)=a-b-1$, then at most one of
$u_2,v_2,w_2$ has degree $a-b-1$.

Otherwise, we may assume that $d(u_2)=d(v_2)=a-b-1$. By Claim 1,
$d(v)=a-b-1$ for any $v\in V(P_1)\cup V(P_2)\backslash\{u_1,u_n\}$.
Then, all vertices on the cycle $C=P_1\cup P_2=u_1u_2\cdots
u_nv_{m-1}\cdots v_2u_1$ have degree $a-b-1$, a contradiction to
Lemma 4. \hfill $\Box$

{\bf Claim 4}. $d(u_1)=d(u_n)=3$.

We only need to prove that $d(u_1)=3$ by Claim 2. If $d(u_1)\neq
3=d_{G_0}(u_1)$, then $d(u_1)=a-b-1$ by Lemma 2, and
\begin{equation}\label{eq:7}
a-b-1>3,\,\,\,\,\mbox{i.e.,}\,\,\, a-b\geq 5.
\end{equation}
From Claim 3, at least two of $u_2,v_2,w_2$ have degree 2. Without
loss of generality, assume that $d(u_2)=d(v_2)=2$. By
Eq.(\ref{eq:1}), $s(u_1)=-(a-b-1)^2+a(a-b-1)-b=(b+1)(a-b-1)-b$. On
the other hand,
$s(u_1)=(d(u_1)-d_{G_0}(u_1))+d(u_2)+d(v_2)+d(w_2)=(a-b)+d(w_2)$ and
$2\leq d(w_2)\leq a-b-1$. We have $(a-b-1)+3\leq (b+1)(a-b-1)-b\leq
2(a-b-1)+1$. And (i) $b+3\leq b(a-b-1)$, implying that $b\neq 0$ and
$b\geq 1$; (ii) $(b-1)(a-b-2)\leq 2$, implying that $b\leq 1$. So,
$b=1$.

By Eq.(\ref{eq:1}), $s(u_2)=-2^2+2a-b=2a-b-4$. And
$d(u_3)=s(u_2)-d(u_1)=a-3<a-b-1$, $d(u_3)=d_{G_0}(u_3)=2$ from Lemma
2 and $u_3\neq u_n$ since $d(u_n)=d(u_1)=a-b-1$. So, $a=d(u_3)+3=5$.
But now $a-b=4$ contradicts Eq.(\ref{eq:7}). Thus, $d(u_1)=3$.
\hfill $\Box$

In the following, we will show that $G$ is not a 2-walk $(a,
b)$-parabolic graph if $G_0=F_3$ by Claims 1-4.

If there is a 2-walk $(a, b)$-parabolic graph $G$ such that
$G_0=F_3$, let $x$ be a pendent vertex of $G$ and $xy$ is the only
edge incident with $x$. Then, $y\in (V(P_1)\cup V(P_2)\cup
V(P_3))\backslash \{u_1,u_n\}$ by Claim 4 and Lemma 3,
$d_{G_0}(y)=2$ and $d(y)=a-b-1$ by Lemma 2.  Let
$N_{G_0}(y)=\{y',y''\}$, $4\leq d_{G_0}(y')+d_{G_0}(y'')\leq
d(y')+d(y'')\leq 2(a-b-1)$ by Lemma 2. By Eq.(\ref{eq:1}),
$s(y)=d(y')+d(y'')+(d(y)-d_{G_0}(y))=-(a-b-1)^2+a(a-b-1)-b=(b+1)(a-b-1)-b$,
i.e., $d(y')+d(y'')=b(a-b-2)+2$. So, $4\leq b(a-b-2)+2\leq
2(a-b-1)$. The left inequality implies $b>0$, and the right
inequality implies that $(b-2)(a-b-2)\leq 0$, i.e., $b\leq 2$, since
$a-b-2>0$ by Lemma 3. So, $b=1$ or $b=2$.

Without loss of generality, we assume $y\in V(P_1)\backslash
\{u_1,u_n\}$. By Claim 1, $d(u_2)=\cdots=d(u_{n-1})=d(y)=a-b-1$.
Using Eq.(\ref{eq:1}),
$$a=\frac{s(u_2)-s(x)}{d(u_2)-d(x)}+d(u_2)+d(x)
=\frac{d(u_3)+1}{a-b-2}+a-b $$ and
\begin{equation}\label{eq:8}
\qquad d(u_3)=b(a-b-2)-1
\end{equation}

(i) If $b=1$, then by Eq.(\ref{eq:8}), $d(u_3)=a-4<a-b-1$, and
$u_3=u_n$, $d(u_3)=d(u_n)=3$. So, $a=7$. By Eq.(\ref{eq:1}), we have
$s(u_1)=d(u_2)+d(v_2)+d(w_2)=-d^2(u_1)+ad(u_1)-b=-9+21-1=11$, and
$d(v_2)+d(w_2)=s(u_1)-d(u_2)=11-5=6$. Since
$d(v_2),d(w_2)\in\{d_{G_0}(v_2),d_{G_0}(w_2), a-b-1\}$ by Lemma 2
and $a-b-1=5$, $d(v_2)=d_{G_0}(v_2)\leq 3$ and
$d(w_2)=d_{G_0}(w_2)\leq 3$. So, $d(v_2)=d_{G_0}(v_2)=3$,
$d(w_2)=d_{G_0}(w_2)=3$, and $v_2=u_n=w_2$, i.e., $m=t=2$,
contradicts that $G$ is a simple graph.

(ii) If $b=2$, then by Eq.(\ref{eq:8}), $d(u_3)=2(a-b-2)-1=2a-9$. If
$n=3$, then $d(u_3)=d(u_n)=3$, and $2a-9=3$ which implies $a=6$.
Otherwise, $d_{G_0}(u_3)=2\neq 2a-9$, and $d(u_3)=a-b-1=a-3$ by
Lemma 2. So, $2a-9=a-3$ and $a=6$. In a word, we have $a=6$, $b=2$,
$d(u_3)=3$. Now, $d(u_2)=d(u_3)=\cdots=d(u_{n-1})=3$ from Claim 1,
$d(v_2)=\cdots=d(v_{m-1})=2$ and $d(w_2)=\cdots=d(w_{t-1})=2$ from
Claims 2-4. Moreover, $m=t=3$ from Lemma 5, i.e., $G\cong H$. Also,
it is easy to check that $H$ is a 2-walk $(6,2)$-parabolic graph.

The proof of Theorem 1 is completed. \hfill $\Box$

{\bf Theorem 2}. A bicyclic graph $G$ with $\delta(G)>1$ has exactly
two main signless Laplacian eigenvalues if and only if $G$ is
isomorphic to one of $G_1,G_2,\cdots,G_7$ (see Figure 4).

{\bf Proof}. Let $G$ be a bicyclic graph. If $G$ is a 2-walk $(a,
b)$-parabolic graph and $\delta(G)>1$, then $G$ is one of the three
graphs $F_1, F_2, F_3$ shown in Figure 1.

{\bf Case 1}. $G=F_1$. Then $l(C_1)=l(C_2)=3$ from Lemma 6(ii),
i.e., $G=G_1$. Using Eq.(\ref{eq:1}), it is easy to check that $G_1$
is a 2-walk $(7,4)$-parabolic graph.

{\bf Case 2}. $G=F_2$. Then $l(C_1)=l(C_2)=3$ from Lemma 6(ii), and
$l(P)=1$ or $l(P)=3$ from Lemma 5(i). So, $G=G_2$ or $G=G_3$. It can
be checked easily that $G_2$ and $G_3$ are 2-walk $(7,5)$-parabolic
and $(6,3)$-parabolic graphs, respectively.

{\bf Case 3}. $G=F_3$. Without loss of generality, we assume that
$l(P_1)\geq l(P_2)\geq l(P_3)$.  Then $l(P_1)\geq l(P_2)\geq 2$
since $G$ is a simple graph.  By Lemma 5(i), $l(P_i)\leq 3$
($i=1,2,3$) and if one of $l(P_1),l(P_2),l(P_3)$ is 3, then the
other two are not 2. So, $G=G_4$ or $G_5$ if $l(P_3)=1$; $G=G_6$ if
$l(P_3)=2$; $G=G_7$ if $l(P_3)=2$. Also, $G_4,G=G_5,G=G_6,G=G_7$ are
2-walk $(6,2)$-parabolic, $(7,5)$-parabolic, $(5,0)$-parabolic and
$(6,3)$-parabolic graphs.

Therefore, a bicyclic graph $G$ with $\delta(G)>1$ has exactly two
main signless Laplacian eigenvalues has exactly two main signless
Laplacian eigenvalues if and only if $G$ is isomorphic to one of
$G_1,G_2,\cdots,G_7$. \hfill $\Box$

\end{document}